\documentclass[conference]{IEEEtran}

\usepackage{cite}
\usepackage{subfigure}
\usepackage{graphicx}
\usepackage{algorithmic}
\usepackage{array}
\usepackage[caption=false,font=footnotesize]{subfig}
\usepackage{url}
\usepackage{amsmath,amssymb,amsfonts}
\usepackage{textcomp}
\usepackage{xcolor,multirow,multicol}

\usepackage{tikztensor}
\usetikzlibrary{calc}
\definecolor{kulblue}{RGB}{0,85,165}
\tikzset{fancy/.style={inner color=kulblue!5!white, outer color=kulblue!20!white, fill faces, opacity=0.95}}
\colorlet{kulblue5}{kulblue!5!white}
\colorlet{kulblue20}{kulblue!20!white}
\colorlet{kulblue30}{kulblue!30!white}
\colorlet{kulblue70}{kulblue!70!white}
\colorlet{kulblue30}{kulblue!30!black}
\colorlet{kulblue60}{kulblue!60!black}
\colorlet{kulblue90}{kulblue!90!black}
\colorlet{myRed}{red!40!white}
    
\newcommand{\bc}[1]{\mbox{\boldmath $\mathcal{#1}$}}
\newcommand{\mb}[1]{\mathbf{#1}}
\newcommand{\F}{\mathrm{F}}
\newcommand{\T}{\mathrm{T}}

\begin{document}

\title{Online Rank-Revealing Block-Term Tensor Decomposition}

\author{\IEEEauthorblockN{Athanasios A. Rontogiannis}
\IEEEauthorblockA{\textit{IAASARS, National Observatory of Athens} \\
152 36 Penteli, Greece \\
tronto@noa.gr}
\and
\IEEEauthorblockN{Eleftherios Kofidis}
\IEEEauthorblockA{\textit{Dept. of Statistics and Insurance Science} \\
\textit{University of Piraeus}\\
185 34 Piraeus, Greece \\
kofidis@unipi.gr}
\and
\IEEEauthorblockN{Paris V. Giampouras}
\IEEEauthorblockA{\textit{Mathematical Inst. for Data Science} \\
\textit{Johns Hopkins University}\\
Baltimore, MD 21218, USA \\
parisg@jhu.edu}

}

\maketitle

\begin{abstract}
The so-called block-term decomposition (BTD) tensor model, especially in its rank-$(L_r,L_r,1)$ version, has been recently receiving increasing attention due to its enhanced ability of representing systems and signals that are composed of \emph{block} components of rank higher than one, a scenario encountered in numerous and diverse applications. Its uniqueness and approximation have thus been thoroughly studied. The challenging problem of estimating the BTD model structure, namely the number of block terms (rank) and their individual (block) ranks, is of crucial importance in practice and has only recently started to attract significant attention. In data-streaming scenarios and/or big data applications, where the tensor dimension in one of its modes grows in time or can only be processed incrementally,  it is essential to be able to perform model selection and computation in a recursive (incremental/online) manner. To date there is only one such work in the literature concerning the (general rank-$(L,M,N)$) BTD model, which proposes an incremental method, however with the BTD rank and block ranks assumed to be a-priori known and time invariant. In this preprint, a novel approach to rank-$(L_r,L_r,1)$ BTD model selection and tracking is proposed, based on the idea of imposing column sparsity jointly on the factors and estimating the ranks as the numbers of factor columns of nonnegligible magnitude. An online method of the alternating iteratively reweighted least squares (IRLS) type is developed and shown to be computationally efficient and fast converging, also allowing the model ranks to change in time. Its time and memory efficiency are evaluated and favorably compared with those of the batch approach. Simulation results are reported that demonstrate the effectiveness of the proposed scheme in both selecting and tracking the correct BTD model. 
\end{abstract}

\section{Introduction}
\label{sec:intro}

Tensors and their decomposition models and methods~\cite{sdfhpf17} have attracted significant attention in numerous application areas in view of their unique ability to represent (explicitly or implicitly) multi-dimensional data and systems and their latent structure. The Tucker decomposition (TD), consisting of a core tensor multiplied by matrix factors on all or some of its modes, and its specialization with super-diagonal core, namely the Canonical Polyadic Decomposition (CPD) (or PARAllel FACtor analysis (PARAFAC or PARAFAC1)), have been the most well-known and studied tensor decomposition models~\cite{sdfhpf17}, in both static and online settings. 
Block-Term Decomposition (BTD) was introduced in~\cite{ldl08b} as a tensor model that combines the CPD and the TD models, in the sense that it decomposes a tensor in a sum of tensors (block terms) that have low multilinear rank (instead of rank one as in CPD). Hence a BTD can be seen as a constrained TD, with its core tensor being block diagonal (see~\cite[Fig.~2.3]{ldl08b}). It can also be seen as a constrained CPD having factors with (some) collinear columns~\cite{ldl08b}. In a way, BTD lies between the two extremes (in terms of core tensor structure), CPD and TD, and it is useful to recall here the related remark made in~\cite{ldl08b}, namely that ``the rank of a higher-order tensor is actually a combination of the two aspects: one should specify the number of blocks \emph{and} their size". 
Accurately and efficiently estimating these numbers and tracking them in time along with the model factors as the tensor grows in one of its modes (common scenario in data-streaming and big data applications) is the aim of this work.

\subsection{Offline BTD model selection and computation}

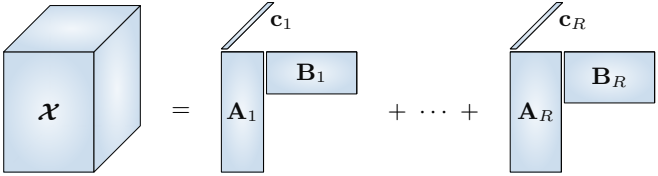
\begin{figure}
\centering
\begin{tikzpicture}[node distance=0.3cm, chain,
term/.style={dim={4,3,4},fancy,tensor scale=0.4}]
\small
% Tensor
\node [tensor, term, dashed back lines] {$\bc{X}$};
% Equality
\node {$=$};
% Term 1
\node [rank-LL1 tensor, term, L=1.4, 2d,
labels={$\mathbf{A}_{1}$}{$\mathbf{B}_{1}$}{$\mathbf{c}_{1}$}] {};
% Plus
\node {$+\ \cdots\ +$};
% Term 2
\node [rank-LL1 tensor, term, L=1.7, 2d,
labels={$\mathbf{A}_{R}$}{$\mathbf{B}_{R}$}{$\mathbf{c}_{R}$}] {};
\end{tikzpicture}
\caption{Rank-$(L_r,L_r,1)$ block-term decomposition.}
\label{fig:BTD}
\end{figure}

Although~\cite{ldl08b} introduced BTD as a sum of $R$ rank-$(L_r,M_r,N_r)$ terms ($r=1,2,\ldots,R$) in general, the special case of rank-$(L_r,L_r,1)$ BTD has attracted a lot more of attention, because of both its more frequent occurrence in a wide range of applications and the existence of more concrete and easier to check uniqueness conditions (cf.~\cite{rkg21} for an extensive review). Note that the mode with the unit modal rank is most commonly associated with the temporal dimension, as it is the case, for example, in tensorial functional Magnetic Resonance Imaging (fMRI)~\cite{ckmt19}. This special yet very popular BTD model is at the focus of the present work and is briefly defined as follows. Consider a 3rd-order tensor, $\bc{X}\in\mathbb{C}^{I\times J\times K}$. Then its rank-$(L_r,L_r,1)$ BTD is written as
\begin{equation}
\bc{X}=\sum_{r=1}^{R}\mb{E}_{r}\circ \mb{c}_{r},
\label{eq:BTD1}
\end{equation}
where $\mb{E}_{r}$ is an $I\times J$ matrix of rank $L_r$, $\mb{c}_{r}$ is a nonzero column $K$-vector and $\circ$ denotes outer product. 
Clearly, $\mb{E}_{r}$ can be written as a matrix product $\mb{A}_{r}\mb{B}_{r}^{\T}$ with the matrices $\mb{A}_{r}\in\mathbb{C}^{I\times L_r}$ and $\mb{B}_{r}\in\mathbb{C}^{J\times L_r}$ being of full column rank, $L_r$. Eq.~(\ref{eq:BTD1}) can thus be re-written as 
\begin{equation}
\bc{X}=\sum_{r=1}^{R}\mb{A}_{r}\mb{B}_{r}^{\T}\circ \mb{c}_{r}.
\label{eq:BTD}
\end{equation}
A schematic representation of the rank-$(L_r,L_r,1)$ BTD is given in Fig.~\ref{fig:BTD}.
It should be apparent from~(\ref{eq:BTD}) and Fig.~\ref{fig:BTD} that CPD results as a special case with all $L_r, r=1,2,\ldots,R$ equal to~1.

In general, $R$ and $L_r$, $r=1,2,\ldots,R$ are assumed \emph{a-priori} known (and it is commonly assumed that all $L_r$ are all equal to $L$, for simplicity). However, unless external information is given (such as in a telecommunications ~\cite{ldl12} or a hyperspectral image (HSI) unmixing application with given or estimated ground truth~\cite{qxzzt17}), there is no way to know these values beforehand. Although overestimation of the block ranks $L_r$s has been observed not to be harmful in some blind source separation applications (e.g., \cite{ldl12}), this is not the case in general~\cite{rkg21}. Besides, in addition to increasing the computational complexity, setting $L_r$ too high may hinder interpretation of the results through letting noise/artifact sources interfere with the desired sources. This holds for $R$ as well, whose choice is known to be more crucial to the obtained performance as it represents the number of ``factors" that generate the data and its over/under-estimation will lead to over/under-fitting, with undesired consequences for the interpretability of the results (cf.~\cite{rkg21} for related references). 

It is known that computing the number of rank-1 terms in a CPD model (i.e. the tensor rank) is NP-hard~\cite{hl13}.
Model selection  for BTD is clearly more challenging than in CPD and TD models and has only recently started to be studied (cf.~\cite{rkg21} for an extensive review of heuristic approaches and techniques). The most recent contribution of this kind can be found in our work~\cite{rkg21}, which relies on a regularization of the squared approximation error function with the sum of the Frobenius norms of the factors reweighted by a diagonal weighting which jointly depends on the factors in two levels: the reweighted norms of $\mb{A}\triangleq \left[\begin{array}{cccc} \mb{A}_1 & \mb{A}_2 & \cdots & \mb{A}_R \end{array}\right]$ and $\mb{B}\triangleq \left[\begin{array}{cccc} \mb{B}_1 & \mb{B}_2 & \cdots & \mb{B}_R \end{array}\right]$ are combined and then coupled with the reweighted norm of $\mb{C}\triangleq \left[\begin{array}{cccc} \mb{c}_1 & \mb{c}_2 & \cdots & \mb{c}_R \end{array}\right]$. This two-level coupling naturally matches the structure of the model in~(\ref{eq:BTD}), making explicit the different roles of $\mb{A},\mb{B}$ and $\mb{C}$. This way, column sparsity is imposed \emph{jointly} on the factors and in a \emph{hierarchical} manner, which allows to estimate the ranks as the numbers of factor columns of non-negligible energy. Following a block coordinate descent (BCD) solution approach, an alternating hierarchical iterative reweighted least squares (HIRLS) algorithm, called BTD-HIRLS, was developed in~\cite{rkg21} that manages to both reveal the ranks and compute the BTD factors at a high convergence rate and low computational cost.

\subsection{Our contribution}

In practice, data may be streaming (arriving sequentially in time) or the data-generation mechanism (the model) may change with time (cf.~\cite{dsbm20,sdplg20} and references therein). In big data applications, the tensor to be decomposed may be too large to fit in the memory and being processed as a whole. A possible workaround is then to recursively update the decomposition model in an incremental manner. In both cases, a sequence of optimization problems results~\cite{dsbm20}. Instead of re-selecting and re-computing the model of the entire tensor every time new data arrive, which is computationally and memory costly even when the previous model is used to initialize the computations and renders the task intractable for fast streamed and/or large-scale tensors, a recursive update is more desirable, namely one that will update the model with relatively few additional operations and will be memory efficient. The aim of such a decomposition `on the fly'~\cite{mmg15} is to track the model in (nearly) real-time with  minimal memory overhead, while attaining a modeling performance comparable to that of the batch (in terms of decomposing the tensor in increasing batches or in its entirety) approach. In any case, an online approach is the only effective way to follow if the model structure (ranks) and parameters change with time as it is the case in numerous application contexts such as video surveillance, network monitoring, dynamic neuroimaging, and others (see the following subsection).

Assume a noisy version of $\bc{X}$ above and that new $I\times J$ frontal slices are being added in its third mode. Assume moreover that the uniqueness conditions hold for all these tensor batches and that $\mb{A},\mb{B}$ (and hence the corresponding multilinear subspace) change only slowly with time. In other words, the assumption is that the incoming data do not have a detrimental effect on the model (i.e., the BTD model keeps being valid, albeit with possibly changing ranks), affecting smoothly only its local features (in the sense of~\cite{zeb18}). Hence (\ref{eq:BTD}) holds throughout, with varying factors $\mb{A},\mb{B},\mb{C}$ and possibly $R$ and $L_r$s. It is natural to assume, for a rank-$(L_r,L_r,1)$ BTD model, that it is the 3rd dimension that grows, while the other dimensions are kept unchanged. The growing mode could correspond to, for example, time in fMRI and dynamic MRI~\cite{mmg15}, spectrum in HSI, etc.
$R$ and $L_r$s are assumed \emph{a-priori unknown} and possibly \emph{time varying}. This is in contrast to all (but one, to be reviewed in the sequel) the existing works on online tensor factorization (OTF).

Our contribution in this work lies in extending the batch approach of~\cite{rkg21} to OTF with possibly changing ranks. Specifically, we propose and develop here a relaxed variant of BTD-HIRLS, which is more amenable to an online extension. The hierarchical nature of BTD-HIRLS is relaxed, considering the third mode separately from the other modes. For extending to the online setting this so-called BTD-IRLS algorithm, we expand on earlier work of ours on online rank-revealing matrix decomposition~\cite{rgk20}. We consider exponential windowing to realize fading memory, that is, gradually forget the effect of past slices. Alternatively, we could employ sliding (truncated) windowing~\cite{ns09}. Whatever is the window type adopted, its duration is to be chosen according to the dynamics of the system under study and can also be adaptive (e.g., \cite{jr19}).
The proposed algorithm, named \emph{Online BTD-IRLS}, is shown to enjoy computational efficiency and fast convergence, inherited from its batch counterpart. Furthermore, it is highly time- and memory-efficient in the sense that the involved quantities are updated recursively in an efficient manner and its computational and memory complexity requirements do not depend on the horizon of the growing dimension. 
Moreover, as demonstrated via simulations, the quality of the data approximation is comparable with that of the batch version and the ranks are estimated and tracked correctly with a high probability. In this presentation, adaptation is done per slice. The method can be easily extended so as to update the model for each new chunk of possibly more than one slices. The algorithm development can be cast in the so-called block successive upper-bound minimization (BSUM) framework~\cite{hrlp16} (see also the appendices) and can be seen as a recursive variant of BTD-IRLS. 

\subsection{Related work}

This subsection presents an extensive overview of the incremental/online tensor decomposition literature with the aim of offering an as complete as possible picture of the framework in which the proposed method is placed. This is not intended to be exhaustive, although most of the related work is included. The interested reader is invited to check the cited works and the references therein. 

Processing a tensor that grows in one of its modes through the streaming of new slices in a way that takes this into account and does not re-compute the model parameters from scratch for every increment in the size of the tensor first appeared in the tensor-based data mining literature~\cite{stpyf08} with the name of Incremental Tensor Analysis (ITA). ITA and its variants (Dynamic Tensor Analysis (DTA), Streaming Tensor Analysis (STA), and Windowed Tensor Analysis (WTA)) were aimed at performing TD on a sequence of tensors and were later extended to incremental higher-order singular value decomposition (HO-SVD) (based on incremental SVD) for the purposes of data mining in intelligent transportation systems~\cite{khyllm14}, computer vision~\cite{lhzzl07,wglt09,hlzsmz11,sbbz14,sjjbz15}, recommendation systems (see also \cite[Chapter~6]{sz16})~\cite{zltc14,zslg14}, handwritten digit recognition~\cite{lyzlllz14}, and epidemics~\cite{kklmnr20} and 
social networks~\cite{ffg20} analysis. Pairwise interactive tensor factorization (PITF), a special case of TD, was studied in this context in~\cite{rvjb19}. A quite special idea for online TD was reported in~\cite{ycl15}. One augments the current TD factors with random columns, orthogonalizes them, finds the TD core tensor with these factors for the newly arrived tensor, reduces the resulting core to the previous dimensions, and combines the factors of the last two decompositions. The idea is that the new decomposition can result from the previous one by updating the factors with random columns and refine them by projecting onto the new data. The use of computationally expensive SVDs is avoided in~\cite{xwmg18} by exploiting the properties of block tensor matrix multiplication to efficiently realize online TD.
In~\cite{kl10}, the tensor is seen as a function of time and the manifold properties of the set of low-rank Tucker decomposable tensors are exploited to formulate the dynamic (through time) Tucker approximation problem with the aid of a system of nonlinear tensor differential equations that are to be numerically solved. Online HO-SVD was also considered in~\cite{vl19} for tracking of the subspaces of the non-evolving modes. Tensor switch, i.e., abrupt changes in the subspaces dimensions, is allowed, with modal rank adjustment based on thresholding the singular values of a rank-revealing URV decomposition. Of course one has to properly pre-set the threshold values, depending on the expected noise level. 

Online CPD was also extensively studied after the pioneering work of~\cite{ns09}; see, e.g., \cite{zvbjd16,nal17,asz20b}. The basic idea (proposed in~\cite{ns09}) is, given a new frontal slice, to find the new row of the $\mb{C}$ factor, use this to find the updated subspace spanned by modes~1 and~2 and finally split it (one way or another) into its $\mb{A}$ and $\mb{B}$ factors. If there is no need to explicitly track $\mb{A}$ and $\mb{B}$, one may simply track the corresponding subspace only. Variations include second-order stochastic gradient (for faster convergence) and updating only one column of the subspace basis matrix at each time (for linear complexity~\cite{nal17,z19}) or use of momentum (Nesterov) stochastic gradient descent (SGD) and artificial noisy perturbations for escaping saddle points~\cite{asz20a}. Online nonlinear least squares (NLS), with possibly dynamic CPD rank, which however needs to be known at every time instant, was developed in~\cite{vvl17} and was shown to outperform~\cite{ns09} in terms of approximation error and speed. The tensor is replaced by its current CPD approximation so that an infinite time horizon can be more easily coped with. All the previous methods are deterministic. A probabilistic online CPD scheme was proposed in~\cite{dzlz18}, which follows a streaming variational Bayesian inference approach to quantify the uncertainty in the CPD model parameters and hence be able to predict missing entries that may arrive in different orders in time, and in all modes.

Constraints, which impose smoothness in time and nonnegativity or sparsity, may also be included depending on the application context, such as, for example, in nonnegative online CPD for time-evolving topic modeling (with possible applications in the analysis of social media-generated data on the Covid-19 pandemic)~\cite{kklmnr20} or in tensor dictionary learning~\cite{tbr19}. Possible solution approaches include alternating optimization with alternating direction method of multipliers (AO-ADMM)~\cite{shsk18}. Incorporating $\ell_1$-norm constraints in dynamic TD allows for rejecting outliers or detecting subspace changes~\cite{cdpm20}. \cite{nhy19} relies on the well-known in robust principal component analysis (PCA)~\cite{xzq19,rgk20} low rank plus sparse representation model to come up with an online CPD scheme for (ADMM-based) outlier-resistant tracking and completion.
For symmetric (moment) tensors, online (SGD) versions of the symmetric tensor power method~\cite{kr02} were developed in~\cite{hnha15,wa16} with applications in latent model learning for community detection and topic modeling.

As mentioned previously, OTF has been also studied in the context of dictionary learning (DL), with the dictionary being Kronecker or Khatri-Rao (KR) product-structured. Thus, an online CPD method based on matrix DL for the tensor unfolding corresponding to the sparse mode was developed in~\cite{rlh20} for 3-way tensors based on the idea that the KR structure of the dictionary can be exploited to disentangle its two factors via least squares (LS) KR factorization (LS-KRF). This was extended in~\cite{sln21} to higher-order tensors through an online CPD scheme that can also incorporate additional constraints and follows an idea analogous to that of~\cite{mbps10} for online matrix DL, which consists of a sparse coding step followed by circular BCD over the factors of the dictionary. BCD with diminishing radius (DR)~\cite{l21} is employed for increased stability. Notably, the development of the scheme of~\cite{sln21} was followed by a convergence proof (also valid for Markovian data), which is rarely the case in OTF works, and was successfully applied in brain video tracking, among other applications. It should be noted that the basic mechanism of~\cite{sln21} (updating of the evolving mode in alternation with BCD over the rest of the modes) is what underlies the scheme proposed in the present work as well and hence the proof in~\cite{sln21} can be extended towards a convergence analysis of our online method and its DL generalization. It should also be noticed that our problem is basically one of sparse factorization with a-priori unknown \emph{support} (i.e., set of the indices $r$ of the nonzero columns of $\mb{C}$ and the indices $l$ of the nonzero columns of the corresponding $\mb{A}_r,\mb{B}_r$ blocks). Thus, the convergence theory of online sparsity learning methods~\cite[Chapter~10]{t20} is also relevant here.

To cope with the large scale of tensors in big data streaming applications, divide-and-conquer approaches have been followed, which work incrementally but in parallel over small batches, fusing the individual CPDs to form that of the overall big tensor. \cite{gpp18} is such an example. CPDs of subtensors that include samples of the newly arrived slices are computed in parallel. Instead of operating on the full data, the method operates on summaries of the data. The rank is estimated with the aid of the Core Consistency Diagnostic (CorConDia) heuristic and a quality control is incorporated to check whether the rank of the new slices disagrees with the current one. The case where the new rank is lower is addressed by solving an assignment problem~\cite{gpp18}. A similar idea, inspired from the philosophy of~\cite{spf14}, underlies the method of~\cite{gpyp19}. 
Similar ideas are found in~\cite{ctzwl21}, where a Bayesian scheme is also developed to perform a statistical data-driven initialization. \cite{hcs16} performs online CPD per blocks, corresponding to those subtensors that are affected by the new data, and accordingly refines the overall CPD. In~\cite{yy20}, the CPD is incremented without having to generate subtensors. \cite{myw18} develops a variant of the method of~\cite{zvbjd16} that employs randomized least squares to more effectively cope with large-scale tensors. A randomization-based online CPD algorithm of the Newton type appears in~\cite{tath21}.

In general, the tensor rank is assumed time invariant and very often \emph{a-priori} known. \cite{pgp18} is a rare exception, where the rank is allowed to change with time (and may correspond to a so-called drift of concepts). The rank of the newly arrived tensor is estimated and an algorithm based on correlating the current CPD factors with those of the new CPD is applied in order to find new concepts or missing concepts. Therefore, a new rank estimation and subsequent CPD procedure is needed for each newly arrived tensor, which considerably complicates the algorithm. Note that in our method rank estimation and tracking is done automatically.  

Besides TD and CPD, alternative tensor decomposition models have been implemented via online methods. These include BTD~\cite{gp20},  PARAFAC2~\cite{gtp20}, and tensor trains~\cite{tatb20}. The so-called \emph{OnlineBTD} algorithm (inspired from the OnlineCPD of~\cite{zvbjd16}) recently reported in~\cite{gp20} is to the best of our knowledge the only method for BTD that operates in an online fashion. It is shown to outperform batch alternating least squares (ALS) and NLS in terms of time and memory efficiency while attaining a comparable approximation error. Moreover, it is made to work for general rank-$(L,M,N)$ BTD and for tensors of higher (than~3) order. However, both the number of block terms and their multilinear ranks are assumed to be fixed and \emph{a-priori} known. The authors show that their method's performance is rather robust to overestimates of these ranks. However, this may not be always the case depending on the application (see~\cite{rgk21} and references therein) plus that one might want to have a sufficiently accurate estimate of the ranks for the purposes of interpreting the data (as in, e.g., HSI, where the rank signifies the number of endmembers and the block ranks stand for the ranks of the corresponding abundance maps). Moreover, it may be the case in practice that the ranks vary with time. Our online method for rank-$(L_r,L_r,1)$ BTD automatically estimates the ranks and tracks them in time. 

In most online methods, the tensor grows in only one (usually having a time meaning) mode. The more general case of the tensor being incremented in more than one or even all of its modes, usually referred to as multi-aspect streaming tensor analysis, can be also relevant in applications (e.g., in recommendation systems with growing time and numbers of movies and users) and is certainly more challenging than the single evolving-mode case. Examples of works that have studied this more general scenario include~\cite{fg15,shgch17,nmgt18,dzlz18,nhy19}. Notably, \cite{nmgt18} also leverages side information (in the form of linear constraints on the Tucker factors). 

Applications of OTF abound. They include unveiling the topology of evolving networks~\cite{sbg17}, spatio-temporal prediction or image in-painting~\cite{tbr19}, multiple-input multiple-output (MIMO) wireless communications~\cite{ns09,yzl20}, brain imaging~\cite{fem15}, monitoring heart-related features from wearable sensors for multi-lead electro-cardiography (ECG)~\cite{xzxx20}, anomaly detection in networks and topic modeling~\cite{stpyf08}, structural health monitoring (in an internet of things (IoT) context)~\cite{asz20b}, online cartography (spectrum map reconstruction in cognitive radio networks)~\cite{jr19}, detection of anomalies in the process of 3D printing~\cite{swzs21}, data traffic monitoring in networks~\cite{stpyf08,mmg15}, cardiac MRI~\cite{mmg15}, stream data compression (e.g., in power distribution systems~\cite{zm20} or in video~\cite{mb18}), and online completion~\cite{mmg15,mnla16,k19}, among others.

\subsection{Outline}

The rest of this preprint is organized as follows. The adopted notation is described in the following subsection. The problem is mathematically stated in Section~\ref{sec:problem}, where useful expressions for the tensor unfoldings and slices are also recalled. A relaxed version of the regularization-based criterion for the batch method of~\cite{rkg21} along with the associated iterative BTD-IRLS procedure is presented in Section~\ref{sec:batch}. This serves as the basis for the development of the online method in Section~\ref{sec:online}. Section~\ref{sec:sims} reports and discusses the simulation results. Conclusions are drawn and future work plans are outlined in Section~\ref{sec:concls}.

\subsection{Notation}
\label{subsec:notation}

Lower- and upper-case bold letters are used to denote vectors and matrices, respectively. Higher-order tensors are denoted by upper-case bold calligraphic letters. For a tensor $\bc{X}$, $\mb{X}_{(n)}$ stands for its mode-$n$ unfolding. $\ast$ stands for the Hadamard product and $\otimes$ for the Kronecker product. The Khatri-Rao product is denoted by $\odot$ in its general (partition-wise) version and by $\odot_{\mathrm{c}}$ in its column-wise version. $\circ$ denotes the outer product. The superscript $^{\T}$ stands for transposition. The Matlab indexing notation is adopted. Thus, for example, $\mb{X}(i,:)$ is the $i$th row of the matrix $\mb{X}$ and $\bc{X}(:,:,k)$ is the $k$th mode-3 (frontal) slice of the tensor $\bc{X}$. The identity matrix of order $N$ and the all ones $M\times N$ matrix are respectively denoted by $\mb{I}_N$ and $\mb{1}_{M\times N}$. $\mb{1}_N$ stands for $\mb{1}_{N\times 1}$. The row vectorization and the trace of a matrix $\mb{X}$ are denoted by $\mathrm{vec}(\mb{X})$ and $\mathrm{tr}(\mb{X})$, respectively. $\nabla_{\mb{X}}$ stands for the gradient operator with respect to (w.r.t) $\mb{X}$. $\mathrm{diag}(\mb{x})$ is the diagonal matrix with the vector $\mb{x}$ on its main diagonal. The block diagonal matrix is denoted by $\mathrm{blockdiag}(\cdot)$. The Euclidean vector norm and the Frobenius matrix and tensor norms are denoted by $\|\cdot\|_{2}$ and $\|\cdot\|_{\F}$, respectively. The mixed 1, 2 ($\ell_{1,2}$) norm of a matrix $\mb{X}=\left[\begin{array}{ccc} \mb{x}_1 & \cdots & \mb{x}_N \end{array}\right]$ is defined as $\sum_{n=1}^{N}\|\mb{x}_i\|_2$. $\mathbb{C}$ is the field of complex numbers.

\section{Problem Statement}
\label{sec:problem}

\begin{figure*}
    \centering
\subfigure[]{\begin{tikzpicture}
    \node [tensor, dim={4,3,1}, fill=red, tensor scale = 0.5](Yk){};
    \node at (Yk.A) [tensor, dim={4,3,4}, fancy, tensor scale=0.5,anchor=E](Y){$\bc{Y}$};
    \node at (Y.label front left) [anchor=east] {$I$};
    \node at (Y.label front top) [anchor=south] {$J$};
    \node at (Y.label bottom right) [anchor=north west] {$k$};
\end{tikzpicture}} \hfill
\subfigure[]{\begin{tikzpicture}[node distance=0cm]
\small
\node [rank-LL1 tensor, dim={4,3,5}, L=1.4, 2d, tensor scale=0.5, fill=red, fill opacity=1](Red1){};
\node at ($(Red1)+(-0.1,-0.1)$) [rank-LL1 tensor, dim={4,3,4}, L=1.4, 2d, tensor scale=0.5, opacity=0.95,
mode 1 inner color = kulblue20,
mode 1 outer color = myRed,
mode 2 inner color = kulblue20,
mode 2 outer color = myRed,
mode 3 inner color = kulblue5,
mode 3 outer color = kulblue20,
labels={$\mathbf{A}_{1}$}{$\mathbf{B}_{1}$}{$\mathbf{c}_{1}$}](Blue1){};
\node at ($(Red1)+(1.5,0.2)$) {$+$};
\end{tikzpicture}
\begin{tikzpicture}[node distance=0cm]
\small
\node [rank-LL1 tensor, dim={4,3,5}, L=2, 2d, tensor scale=0.5, fill=red, fill opacity=1](Red2){};
\node at ($(Red1)+(-0.1,-0.1)$) [rank-LL1 tensor, dim={4,3,4}, L=2, 2d, tensor scale=0.5, opacity=0.95,
mode 1 inner color = kulblue20,
mode 1 outer color = myRed,
mode 2 inner color = kulblue20,
mode 2 outer color = myRed,
mode 3 inner color = kulblue5,
mode 3 outer color = kulblue20,
labels={$\mathbf{A}_{2}$}{$\mathbf{B}_{2}$}{$\mathbf{c}_{2}$}](Blue2){};
\node at ($(Red2)+(2.1,0.2)$) {$+\ \cdots\ +$};
\end{tikzpicture}
\begin{tikzpicture}[node distance=0cm]
\small
\node [rank-LL1 tensor, dim={4,3,5}, L=1.7, 2d, tensor scale=0.5, fill=red, fill opacity=1](RedR){};
\node at ($(RedR)+(-0.1,-0.1)$) [rank-LL1 tensor, dim={4,3,4}, L=1.7, 2d, tensor scale=0.5, opacity=0.95,
mode 1 inner color = kulblue20,
mode 1 outer color = myRed,
mode 2 inner color = kulblue20,
mode 2 outer color = myRed,
mode 3 inner color = kulblue5,
mode 3 outer color = kulblue20,
labels={$\mathbf{A}_{R}$}{$\mathbf{B}_{R}$}{$\mathbf{c}_{R}$}](BlueR){};
\end{tikzpicture}
}
\caption{(a) Growth of a 3rd-order tensor in its (evolving) third mode 
with an extra slice (red) and (b) updating of its rank-$(L_r,L_r,1)$ BTD by adding a new row (red) to the third-mode factor matrix
and appropriately modifying the factor matrices in the other modes (pink).}
    \label{fig:incremental}
\end{figure*}
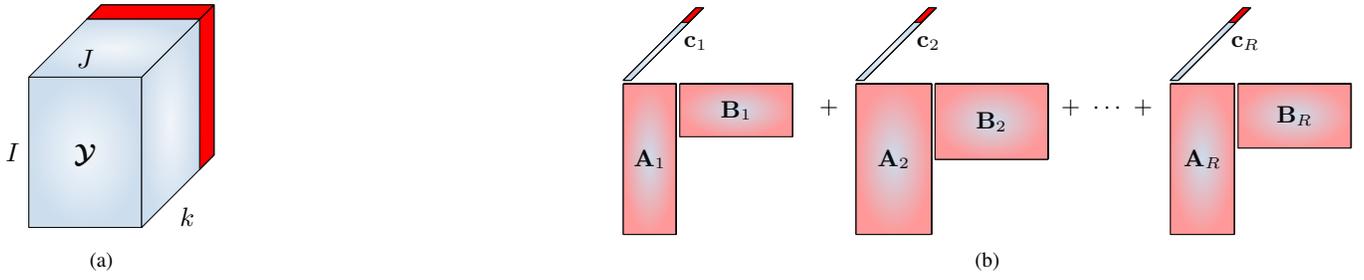

Let the $I\times J\times k$ tensor $\bc{Y}^{(k)}$ grow in its 3rd mode, i.e. for increasing $k$. That is, an additional $I\times J$ frontal slice is considered (in an incremental or streaming fashion) per step  (see Fig.~\ref{fig:incremental}). 
The aim is to compute the best (in the least squares sense) rank-$(L_r,L_r,1)$ approximation of $\bc{Y}^{(k)}$, 
\[
\bc{\hat{Y}}^{(k)}=\sum_{r=1}^R\mb{A}^{(k)}_{r}\mb{B}^{(k)\T}_{r}\circ\mb{c}^{(k)}_{r},
\]
in a recursive manner, that is, based on the BTD model for the $I\times J\times (k-1)$ tensor $\bc{Y}^{(k-1)}$, compute  
the matrices $\mb{A}^{(k)}_{r}=\left[\begin{array}{cccc} \mb{a}^{(k)}_{r,1} & \mb{a}^{(k)}_{r,2} & \cdots & \mb{a}^{(k)}_{r,L_r} \end{array}\right]\in\mathbb{C}^{I\times L_r}$, 
$\mb{B}^{(k)}_r=\left[\begin{array}{cccc} \mb{b}^{(k)}_{r,1} & \mb{b}^{(k)}_{r,2} & \cdots & \mb{b}^{(k)}_{r,L_r} \end{array}\right]\in\mathbb{C}^{J\times L_r}$, and $\mb{C}^{(k)}\in\mathbb{C}^{k\times R}$, with the ranks $R$ and $L_r$, $r=1,2,\ldots,R$ assumed \emph{a-priori unknown and possibly varying with $k$}. 

Recall that, in terms of its mode unfoldings, the tensor in~(\ref{eq:BTD}) can be written as~\cite{ldl08b}
\begin{eqnarray}
\mb{X}_{(1)}^{\T} \!\!\!\!\! & = & \!\!\!\!\! (\mb{B}\odot\mb{C})\mb{A}^{\T}\triangleq \mb{P}\mb{A}^{\T}, \label{eq:X1} \\
\mb{X}_{(2)}^{\T} \!\!\!\!\! & = & \!\!\!\!\! (\mb{C}\odot\mb{A})\mb{B}^{\T}\triangleq \mb{Q}\mb{B}^{\T}, \label{eq:X2} \\
\mb{X}_{(3)}^{\T} \!\!\!\!\! & = & \!\!\!\!\!
\left[\begin{array}{ccc} (\mb{A}_1\odot_{\mathrm{c}} \mb{B}_1)\mb{1}_{L_1} & \cdots & (\mb{A}_R\odot_{\mathrm{c}} \mb{B}_R)\mb{1}_{L_R}\end{array}\right]\mb{C}^{\T} \nonumber \\
\!\!\!\!\!\!\! & \triangleq & \mb{S}\mb{C}^{\T}.
\label{eq:X3}
\end{eqnarray}
Moreover, its $k$th frontal slice can be expressed as
\begin{equation}
\bc{X}(:,:,k)=\mb{A}\mathrm{blockdiag}(c_{k,1}\mb{I}_{L_1},c_{k,2}\mb{I}_{L_2},\ldots,c_{k,R}\mb{I}_{L_R})\mb{B}^{\T},
    \label{eq:fslice1}
\end{equation}
which, for the most common case of all equal $L_r=L$, becomes
\begin{equation}
\bc{X}(:,:,k) = \mb{A}(\mathrm{diag}(\mb{C}(k,:))\otimes \mb{I}_L)\mb{B}^{\T}.
    \label{eq:fslice2}
\end{equation}
Note that, for $L=1$, the above yields the well-known expression for the frontal slices of a CPD-modeled tensor~\cite{sdfhpf17}.
These expressions will be used in the sequel to solve for $\mb{A},\mb{B},\mb{C}$ in a BCD manner. 

\section{Batch BTD-IRLS Algorithm}
\label{sec:batch}

%Setting $\xi=1$ in (\ref{eq:ewminp}), we arrive at the following {\it batch} minimization problem 
The problem solved by BTD-HIRLS in~\cite{rkg21} can be formulated as
\begin{align}
&\underset{\mb{A},\mb{B},\mb{C}}{\mathrm{min}}\frac{1}{2}\left\|\bc{Y} - \sum_{r=1}^R\mb{A}_r\mb{B}_r^{\T}\circ\mb{c}_r\right\|_{\F}^2 + \nonumber \\
 &\lambda\sum_{r=1}^{R}\sqrt{\sum_{l=1}^{L}\sqrt{\|\mb{a}_{r,l}\|_{2}^2 + \|\mb{b}_{r,l}\|_{2}^2 + \eta^2} + \|\mb{c}_r\|_{2}^2 + \eta^2},
\label{eq:BTD-HIRLS} 
\end{align}
where $\eta^2$ is a very small positive constant that ensures smoothness at zero and $R$ and $L$ here stand for the initial {\it (over)estimates} of the model rank parameters. Observe that the definition of the regularizer fully matches the structure of the BTD model in~(\ref{eq:BTD}). For each block term, $r$, the blocks $\mb{A}_r,\mb{B}_r$ are coupled together and, at a higher level, coupled with the corresponding column of $\mb{C}$. Nonetheless, a simpler definition, which relaxes the coupling between $\mb{A}_r,\mb{B}_r$ and $\mb{c}_r$ and greatly facilitates the development of a recursive decomposition scheme, is also possible and is given below:
\begin{align}
&\min_{\mb{A},\mb{B},\mb{C}}\frac{1}{2}\left\|\bc{Y}-\sum_{r=1}^{R}\left(\mb{A}_r\mb{B}_r^{\T}\right)\circ \mb{c}_r\right\|_{\F}^2 + \nonumber \\
&\lambda\sum_{r=1}^{R}\sum_{l=1}^{L}\sqrt{\|\mb{a}_{r,l}\|_{2}^2 + \|\mb{b}_{r,l}\|_{2}^2 + \eta^2} + \mu\sum_{r=1}^{R}\sqrt{\|\mb{c}_r\|_{2}^2 + \eta^2}, 
\label{eq:bminp}
\end{align}
where the regularization parameters of the terms associated with $\mb{A},\mb{B}$ and $\mb{C}$, namely $\lambda$ and $\mu$, may in general differ. 
%The objective function in (\ref{eq:bminp}) is a simplified variant of that minimized in \cite{rkg21} and results by decoupling factor $\mb{C}$  from factors $\mb{A}$ and $\mb{B}$ in the regularizer. 
This modification makes the resulting scheme more flexible and amenable to incremental processing (see Section~\ref{sec:online}). It should also be noted that, ignoring $\eta^2$, the two regularization terms above coincide with the $\ell_{1,2}$ norms of $\left[ \mb{A}^{\T}\;\mb{B}^{\T}\right]^{\T}$ and $\mb{C}$, respectively. Hence, as detailed and demonstrated in~\cite{rkg21}, it is expected that this regularization scheme will promote column sparsity simultaneously on the factors $\mb{A},\mb{B}$ (jointly) and $\mb{C}$, allowing the actual ranks $R$ and $L_r$s to be recovered. Following a similar methodology as in~\cite{rkg21}, namely employing a BCD scheme (with the blocks being the BTD factors $\mb{A},\mb{B},\mb{C}$) and relying on majorization-minimization (MM)~\cite{hl04} for each block, we can develop an iterative reweighted least squares (IRLS)-type algorithm for solving~(\ref{eq:bminp}). The algorithm, called here BTD-IRLS, is tabulated as Algorithm~1 and is seen to only involve closed-form matrix-wise operations for $\mb{A},\mb{B},\mb{C}$ at each iteration.
\begin{table}
\centering
\title{Algorithm 1: The BTD-IRLS algorithm}
 \begin{tabular}{|l|}
 \hline \\
  {\bf Algorithm 1}: BTD-IRLS algorithm\\ \hline 
  Input: $\bc{Y}$,$\lambda,\mu,R_{\mathrm{ini}},L_{\mathrm{ini}}$ \\
  Output: Best (in the sense of~(\ref{eq:bminp})) BTD approximation of $\bc{Y}$ \\
  Initialize: $\mb{A}^{(0)},\mb{B}^{(0)}, \mb{C}^{(0)}$    \\
  $n\leftarrow 0$ \\
  \bf{repeat}\\
	  \hspace{0.3cm} Compute $\mb{D}_1^{(n)},\mb{D}_2^{(n)}$ from~(\ref{eq:Dn1}) and~(\ref{eq:Dn2}) \\
	   \hspace{0.3cm} $\mb{P}^{(n)} \leftarrow \mb{B}^{(n)}\odot\mb{C}^{(n)}$\\
    \hspace{0.3cm} $\mb{A}^{(n+1)} \leftarrow \mb{Y}_{(1)}\mb{P}^{(n)}\left(\mb{P}^{(n){\T}}\mb{P}^{(n)}+\lambda\mb{D}_2^{(n)}\right)^{-1}$ \\
    \hspace{0.3cm} $\mb{Q}^{(n)} \leftarrow \mb{C}^{(n)}\odot\mb{A}^{(n)}$\\
    \hspace{0.3cm} $\mb{B}^{(n+1)} \leftarrow \mb{Y}_{(2)}\mb{Q}^{(n)}\left(\mb{Q}^{(n){\T}}\mb{Q}^{(n)}+\lambda\mb{D}_2^{(n)}\right)^{-1}$ \\
    \hspace{0.3cm} $\mb{S}^{(n)} \leftarrow \left[\begin{array}{ccc} \left(\mb{A}_1^{(n)}\odot_{\mathrm{c}} \mb{B}_1^{(n)}\right)\mb{1}_L & \cdots & \left(\mb{A}_R^{(n)}\odot_{\mathrm{c}} \mb{B}_R^{(n)}\right)\mb{1}_L\end{array}\right]$ \\
    \hspace{0.3cm} $\mb{C}^{(n+1)} \leftarrow \mb{Y}_{(3)}\mb{S}^{(n)}(\mb{S}^{(n){\T}}\mb{S}^{(n)}+\mu\mb{D}_1^{(n)})^{-1}$ \\
    \hspace{0.4cm}$n \leftarrow n+1$\\
    \bf{until} {\it convergence} \\
    \hline
 \end{tabular}
\end{table}
The ranks are over-estimated as $R=R_{\mathrm{ini}}$ and $L_r=L_{\mathrm{ini}}, r=1,2,\ldots,R$ and, provided $\lambda,\mu$ are appropriately selected, their true values are recovered as the numbers of columns of non-negligible magnitude of $\mb{C}$ and the corresponding $\mb{A}_r,\mb{B}_r$ blocks, respectively. This can be done either after convergence or in the course of the procedure, accompanied by the respective column pruning. 
BTD-IRLS differs from the BTD-HIRLS scheme of~\cite{rkg21} in that reweighting is done separately for $\mb{A},\mb{B}$ and $\mb{C}$ (reflecting the `de-coupled' nature of the regularizer) and not in a two-level hierarchy as in BTD-HIRLS, hence the name of the new algorithm. The reweighting diagonal matrices $\mb{D}_1^{(n)}$ and $\mb{D}_2^{(n)}$, of order $R$ and $LR$, respectively, are given by
\begin{equation}
    \mb{D}_1^{(n)}(r,r)=\left(\|\mb{c}_{r}^{(n)}\|_{2}^2 + \eta^2\right)^{-1/2}
    \label{eq:Dn1}
\end{equation}
    and
\begin{align}
    &\mb{D}_2^{(n)}\left((r-1)L+l,(r-1)L+l\right)= \nonumber \\ & \left(\|\mb{a}_{r,l}^{(n)}\|_{2}^2 + \|\mb{b}_{r,l}^{(n)}\|_{2}^2 + \eta^2\right)^{-1/2}
    \label{eq:Dn2}
\end{align}
and are applied for reweighting $\mb{C}$ and $\mb{A},\mb{B}$, respectively.
This modification of BTD-HIRLS inherits its computational efficiency, as it will be detailed in Section~\ref{sec:complexity}. Another difference with the method of~\cite{rkg21} comes from the fact that the two, now `de-coupled', regularization terms can be weighed by generally different regularization parameters, offering the possibility of penalizing the overestimation of the number of block terms and the number of components in each differently. As we will demonstrate via simulations, the above algorithm shares the rank-revelation ability of its counterpart of~\cite{rkg21} and is also fast converging. An online version of it is developed next.

\section{A Rank-Revealing Online BTD Algorithm}
\label{sec:online}

%Exponential windowing, with the forgetting factor $0<\xi\leq 1$, will be assumed. The problem must be solved in an incremental manner, that is, by updating its solution for $\bc{Y}^{(k-1)}$ while minimizing the additional computational and memory cost. In this regard, 
In order to incorporate in the previous method the ability to track time-varying BTD models, we define an exponentially windowed version of the objective function in~(\ref{eq:bminp}). The optimization problem at time step $k$ is formulated as
\begin{align}
&\min_{\mb{A},\mb{B},\mb{C}}\frac{1}{2}\sum_{\kappa = 1}^{k}\xi^{k-\kappa}\left\|\mb{Y}^{(\kappa)}-\mb{A}\left( \mathrm{diag}(\boldsymbol{\gamma}^{(\kappa)})\otimes\mb{I}_{L}\right)\mb{B}^{\T}\right\|_{\F}^2 + \nonumber \\
&\!\!\!\!\!\lambda\sum_{r=1}^{R}\sum_{l=1}^{L}\sqrt{\|\mb{a}_{r,l}\|_{2}^2 + \|\mb{b}_{r,l}\|_{2}^2 + \eta^2}\! +\! \mu\sum_{r=1}^{R}\sqrt{\|\mb{\Xi}^{(k)\frac{1}{2}}\mb{c}_r\|_{2}^2 + \eta^2},
\label{eq:ewminp}
\end{align}
where $\mb{Y}^{(\kappa)}$ is the $\kappa$th $I\times J$ slice of $\bc{Y}^{(k)}$, $\boldsymbol{\gamma}^{(\kappa)\T}\triangleq \mb{C}(\kappa,:)$ is the $\kappa$th row of $\mb{C}$, and $\mb{\Xi}^{(k)} \triangleq \text{diag}\left(\xi^{k-1},\ldots,\xi,1\right)$, with $0<\xi\leq 1$ being the forgetting factor. The aim is, given the model for the tensor consisting of the first $k-1$ slices and the new, $k$th slice $\mb{Y}^{(k)}$, to update the model parameters $\mb{A},\mb{B},\mb{C}$ and the model orders $R$ and $L_r,r=1,2,\ldots,R$ based on their values at $\kappa=k-1$, in a time- and memory-efficient way. As previously, $R$ and $L$ are the over-estimates of the rank parameters, with all $L_r$ being over-estimated as $L$, and we have thus employed the corresponding expression for the frontal slice, eq.~(\ref{eq:fslice2}). It must be emphasized that these parameters are not only a-priori unknown but also potentially changing with $k$.
The first term in~(\ref{eq:ewminp}) is the data fidelity cost, represented by the sum of the exponentially weighted squared errors between the slices of the data tensor and those of its BTD approximation. The rest of the objective is inspired from the regularization part of~(\ref{eq:bminp}). Exponential time-weighting is also applied on the columns of $\mb{C}$, while no weighting is required in the term involving the factors $\mb{A}$ and $\mb{B}$ since these do not change in size with $k$. As it is common in the OTF literature, we will make the assumption that the latter factors, that is those corresponding to the non-evolving modes, are only changing slowly. 

Let $\mb{Y}_{(1)}^{(k-1)},\mb{Y}_{(2)}^{(k-1)}$ and $\mb{Y}_{(3)}^{(k-1)}$ be the mode unfoldings of the $I\times J\times (k-1)$ tensor $\bc{Y}^{(k-1)}$, available at step $k-1$, and  $\mb{y}^{(k)}\triangleq\mathrm{vec}(\mb{Y}^{(k)})$, where $\mb{Y}^{(k)}$ is the new $I\times J$ slice that is included at step $k$. Then (cf.~\cite[Fig.~2]{hlzsmz11}) the mode unfoldings of the incremented tensor $\bc{Y}^{(k)}$ can be expressed as follows:
\begin{eqnarray}
\mb{Y}_{(1)}^{(k)} &=& \left[\begin{array}{cc} \mb{Y}_{(1)}^{(k-1)} & \mb{Y}^{(k)}\end{array}\right]\mb{U}^{(k)},\label{eq:Yk1} \\
\mb{Y}_{(2)}^{(k)} &=& \left[\begin{array}{cc} \mb{Y}_{(2)}^{(k-1)}\; \mb{Y}^{(k)\T} \end{array}\right], \label{eq:Yk2}\\
\mb{Y}_{(3)}^{(k)} &=& \begin{bmatrix} \mb{Y}_{(3)}^{(k-1)} \\ \mb{y}^{(k)\T} \end{bmatrix}, \label{eq:Yk3}
\end{eqnarray}
where $\mb{U}^{(k)}$ is the $kJ\times kJ$ permutation matrix that moves the $j$th column of $\mb{Y}^{(k)}$ to the $jk$th position of the resulting matrix. 

It follows from~(\ref{eq:X3}) and~(\ref{eq:Yk3}) and the assumption of slowly-varying $\mb{A},\mb{B}$ that the factor $\mb{C}$ will only change in its new row at each step $k$, that is,
\begin{eqnarray}
\mb{C}^{(k)} = \begin{bmatrix} \mb{C}^{(k-1)} \\ \boldsymbol{\gamma}^{(k)\T} \end{bmatrix}.
\label{eq:Ck}
\end{eqnarray}
Hence the problem of estimating $\mb{C}$ at the $k$th step can be cast from~(\ref{eq:ewminp}) in terms of its last row and using~(\ref{eq:X3}) can be written as:
\begin{align}
\boldsymbol{\gamma}^{(k)} = \arg\min_{\boldsymbol{\boldsymbol{\gamma}}}&\frac{1}{2}\left\|\mb{y}^{(k)}-\mb{S}^{(k-1)}\boldsymbol{\gamma}\right\|_{2}^2 + \nonumber \\
&\mu\sum_{r=1}^{R}\sqrt{\xi\left\|\mb{\Xi}^{(k-1)\frac{1}{2}}\mb{c}_{r}^{(k-1)}\right\|_{2}^2 + \gamma_r^2 + \eta^2},
\label{eq:objCkp1}    
\end{align}
where $\mb{c}_{r}^{(k-1)}$ is the $r$th column of $\mb{C}^{(k-1)}$ and  $\gamma_r$ is the $r$th element of the optimization  variable  $\boldsymbol{\gamma}$. This sub-problem can be solved via MM (see Appendix~\ref{sec:proofs}), leading to the following closed-form solution for $\boldsymbol{\gamma}^{(k)}$:
\begin{eqnarray}
\boldsymbol{\gamma}^{(k)} = \left(\mb{S}^{(k-1)\T}\mb{S}^{(k-1)}+\mu\mb{D}_{1}^{(k-1)}\right)^{-1}\mb{S}^{(k-1)\T}\mb{y}^{(k)},
\label{eq:Ck+1}
\end{eqnarray}
where (cf.~(\ref{eq:X3})) 
\begin{eqnarray*}
\lefteqn{\mb{S}^{(k-1)}\triangleq} \\
& & \!\!\!\!\!\! 
\left[\begin{array}{ccc} (\mb{A}^{(k-1)}_1\odot_{\mathrm{c}} \mb{B}^{(k-1)}_1)\mb{1}_{L_1} & \cdots & (\mb{A}^{(k-1)}_R\odot_{\mathrm{c}} \mb{B}^{(k-1)}_R)\mb{1}_{L_R}\end{array}\right]
\end{eqnarray*}
and $\mb{D}_1^{(k-1)}$ is the $R\times R$ diagonal matrix with diagonal entries
\begin{equation}
    \mb{D}_1^{(k-1)}(r,r) = \left( \xi\left\|\mb{\Xi}^{(k-1)\frac{1}{2}}\mb{c}_{r}^{(k-1)}\right\|_{2}^2 + \eta^2\right)^{-1/2}. 
    \label{eq:D1}
\end{equation}
Observe that the squared norm above can be recursively computed:
\begin{equation}
\mb{D}_1^{(k)}(r,r)=\left( \xi^2\left\|\mb{\Xi}^{(k-1)\frac{1}{2}}\mb{c}_{r}^{(k-1)}\right\|_{2}^2 + \xi(\gamma_r^{(k)})^2+\eta^2\right)^{-1/2}. 
\label{eq:D1k}
\end{equation}

$\mb{C}$ is updated row-by-row from~(\ref{eq:Ck+1}). On the contrary, each of the factors $\mb{A}$ and $\mb{B}$ must be updated as a whole upon the arrival of a new slice (cf.~Fig.~\ref{fig:incremental}(b)). Writing~(\ref{eq:ewminp}) in terms of the mode-1 unfoldings and making use of~(\ref{eq:X1}) we obtain the $\mb{A}$ sub-problem in the form
\begin{align}
\mb{A}^{(k)} = \arg\min_{\mb{A}}&\frac{1}{2}\left\|{}_{J}\mb{\Xi}^{(k)\frac{1}{2}}\left(\mb{Y}_{(1)}^{(k)\T}-\mb{P}^{(k)}\mb{A}^{\T}\right)\right\|_{\F}^2 + \nonumber \\
&\lambda\sum_{r=1}^{R}\sum_{l=1}^{L}\sqrt{\|\mb{a}_{r,l}\|_{2}^2 + \|\mb{b}_{r,l}^{(k-1)}\|_{2}^2 + \eta^2},
\label{eq:objA}    
\end{align}
where $\mb{P}^{(k)}=\mb{B}^{(k)}\odot\mb{C}^{(k)}$ (cf.~(\ref{eq:X1})) and ${}_{J}\mb{\Xi}^{(k)}$ stands for $\mb{I}_J\otimes\mb{\Xi}^{(k)}$. Solving~(\ref{eq:objA}) via MM (see~Appendix~\ref{sec:proofs}) leads to the following closed-form expression for $\mb{A}^{(k)}$:
\begin{eqnarray}
\lefteqn{\mb{A}^{(k)} =} \nonumber \\ & & \!\!\!\!\!\!\!\!
\mb{Y}_{(1)}^{(k)}{}_{J}\mb{\Xi}^{(k)}\mb{P}^{(k)}\left(\mb{P}^{(k)\T}{}_{J}\mb{\Xi}^{(k)}\mb{P}^{(k)}+\lambda\mb{D}_{2}^{(k-1)}\right)^{-1}, 
\label{eq:Akp1}
\end{eqnarray}
where $\mb{D}_{2}^{(k-1)}$ is the $LR\times LR$ diagonal matrix whose $((r-1)L+l)$th diagonal entry is given by 
\begin{eqnarray}
    &\mb{D}_2^{(k-1)}\left((r-1)L+l,(r-1)L+l\right)  = \nonumber \\ &\left(\|\mb{a}_{r,l}^{(k-1)}\|_{2}^2 + \|\mb{b}_{r,l}^{(k-1)}\|_{2}^2 + \eta^2\right)^{-1/2}. 
    \label{eq:D2}
\end{eqnarray}
With $\mb{V}_{\mb{A}}^{(k)} \triangleq \mb{P}^{(k)\T}{}_{J}\mb{\Xi}^{(k)}\mb{P}^{(k)}$ and $\mb{G}_{\mb{A}}^{(k)} \triangleq \mb{Y}_{(1)}^{(k)}{}_{J}\mb{\Xi}^{(k)}\mb{P}^{(k)}$, (\ref{eq:Akp1}) is more compactly written as 
\begin{eqnarray}
\mb{A}^{(k)} = \mb{G}_{\mb{A}}^{(k)}\left(\mb{V}_{\mb{A}}^{(k)}+\lambda\mb{D}_{2}^{(k-1)}\right)^{-1}. 
\label{eq:Ak}
\end{eqnarray}
In addition, making use of the assumption that $\mb{B}^{(k)} \approx \mb{B}^{(k-1)}$ and recalling the definition of the permutation matrix $\mb{U}^{(k)}$, $\mb{P}^{(k)}$ can be approximately written as
\begin{eqnarray}
\mb{P}^{(k)} &\approx& \mb{U}^{(k)\T}\begin{bmatrix}\mb{B}^{(k-1)}\odot \mb{C}^{(k-1)}  \\ \mb{B}^{(k-1)}\odot \boldsymbol{\gamma}^{(k)\T} \end{bmatrix} \nonumber \\
&= & \mb{U}^{(k)\T}\begin{bmatrix} \mb{P}^{(k-1)}  \\ \mb{B}^{(k-1)}\odot \boldsymbol{\gamma}^{(k)\T} \end{bmatrix}.
\label{eq:Pkp1}
\end{eqnarray}
Since $\boldsymbol{\gamma}^{(k)}$ is a vector, the latter KR product can be equivalently written as 
\begin{equation}
    \mb{B}^{(k-1)}\odot \boldsymbol{\gamma}^{(k)\T}=\mb{B}^{(k-1)}(\mathrm{diag}(\boldsymbol{\gamma}^{(k)})\otimes \mb{I}_L).
    \label{eq:Bgamma}
\end{equation}
Using the following recursion for $\mb{\Xi}$,
\[
\mb{\Xi}^{(k)}=\left[\begin{array}{cc} \xi\mb{\Xi}^{(k-1)} & \mb{0} \\ \mb{0} & 1 \end{array}\right],
\]
and the definition of $\mb{U}^{(k)}$ readily leads to the following
\begin{equation}
\mb{U}^{(k)\T}{}_{J}\mb{\Xi}^{(k)}\mb{U}^{(k)}=\left[\begin{array}{cc} \xi\cdot{}_{J}\mb{\Xi}^{(k-1)} & \mb{0} \\ \mb{0} & \mb{I}_J \end{array}\right].
\label{eq:UXiU}
\end{equation}
From~(\ref{eq:Pkp1}) and~(\ref{eq:Yk1}) and making use of~(\ref{eq:UXiU}) we can easily arrive at the following recursive formulas for updating $\mb{V}_{\mb{A}}^{(k)}$ and $\mb{G}_{\mb{A}}^{(k)}$:
\begin{eqnarray}
\lefteqn{\mb{V}_{\mb{A}}^{(k)} =} \label{eq:HAk} \\
& &  \xi\mb{V}_{\mb{A}}^{(k-1)} + \left(\mb{B}^{(k-1)}\odot\boldsymbol{\gamma}^{(k)\T}\right)^{\T}\left(\mb{B}^{(k-1)}\odot\boldsymbol{\gamma}^{(k)\T}\right) \nonumber
\end{eqnarray}
and
\begin{equation}
\mb{G}_{\mb{A}}^{(k)} = \xi\mb{G}_{\mb{A}}^{(k-1)} + \mb{Y}^{(k)}\left( \mb{B}^{(k-1)}\odot\boldsymbol{\gamma}^{(k)\T}\right).
\label{eq:GAk}
\end{equation}

$\mb{B}$ can be updated at step $k$ by solving the corresponding sub-problem of~(\ref{eq:ewminp}), which, in terms of the mode-2 unfoldings can be expressed as follows:
\begin{align}
\mb{B}^{(k)} = \arg\min_{\mb{B}}&\frac{1}{2}\left\|\mb{\Xi}_{I}^{(k)\frac{1}{2}}\left(\mb{Y}_{(2)}^{(k)\T}-\mb{Q}^{(k)}\mb{B}^{\T}\right)\right\|_{\F}^2 + \nonumber \\
&\lambda\sum_{r=1}^{R}\sum_{l=1}^{L}\sqrt{\|\mb{a}_{r,l}^{(k-1)}\|_{2}^2 + \|\mb{b}_{r,l}\|_{2}^2 + \eta^2},
\label{eq:objB}    
\end{align}
where (cf.~(\ref{eq:X2})) $\mb{Q}^{(k)}\triangleq \mb{C}^{(k)}\odot\mb{A}^{(k)}$ and $\mb{\Xi}_{I}^{(k)} \triangleq 
\mb{\Xi}^{(k)}\otimes\mb{I}_I$. In a manner analogous to the $\mb{A}$ sub-problem, the unique solution for the updated $\mb{B}$ results as
\begin{eqnarray}
\mb{B}^{(k)} = \mb{Y}_{(2)}^{(k)}\mb{\Xi}_{I}^{(k)}\mb{Q}^{(k)}\left(\mb{Q}^{(k)\T}\mb{\Xi}_{I}^{(k)}\mb{Q}^{(k)}+\lambda\mb{D}_{2}^{(k-1)}\right)^{-1}. 
\label{eq:Bkp1}
\end{eqnarray}
Again, employing the assumption that $\mb{A}^{(k)}\approx \mb{A}^{(k-1)}$, $\mb{Q}^{(k)}$ can be approximated as described below:  \begin{eqnarray}
\mb{Q}^{(k)} &\approx&  \begin{bmatrix}\mb{C}^{(k-1)}\odot\mb{A}^{(k-1)}  \\ \boldsymbol{\gamma}^{(k)\T}\odot\mb{A}^{(k)} \end{bmatrix} \nonumber \\
&=&  \begin{bmatrix}\mb{Q}^{(k-1)} \\  \boldsymbol{\gamma}^{(k)\T}\odot\mb{A}^{(k)} \end{bmatrix},
\label{eq:Qkp1}
\end{eqnarray}
where $\boldsymbol{\gamma}^{(k)\T}\odot\mb{A}^{(k)}=\mb{A}^{(k)}
(\mathrm{diag}(\boldsymbol{\gamma}^{(k)})\otimes \mb{I}_L)$.
From~(\ref{eq:Qkp1}) and~(\ref{eq:Yk2}), and using the recursion 
\[
\mb{\Xi}_{I}^{(k)}=\left[\begin{array}{cc} \xi\mb{\Xi}_{I}^{(k-1)} & \mb{0} \\ \mb{0} & \mb{I}_I \end{array}\right],
\]
$\mb{V}_{\mb{B}}^{(k)} \triangleq \mb{Q}^{(k)\T}\mb{\Xi}_{I}^{(k)}\mb{Q}^{(k)}$ and $\mb{G}_{\mb{B}}^{(k)} \triangleq \mb{Y}_{(2)}^{(k)}\mb{\Xi}_{I}^{(k)}\mb{Q}^{(k)}$ can be recursively computed as 
\begin{eqnarray}
\mb{V}_{\mb{B}}^{(k)} \!\!\!\!\!\! & = & \!\!\!\!\!\! \xi\mb{V}_{\mb{B}}^{(k-1)} + \left(\boldsymbol{\gamma}^{(k)\T}\odot\mb{A}^{(k)}\right)^{\T}\left(\boldsymbol{\gamma}^{(k)\T}\odot\mb{A}^{(k)}\right), 
\label{eq:HBk} \\
\mb{G}_{\mb{B}}^{(k)} \!\!\!\!\!\! & = & \!\!\!\!\!\! \xi\mb{G}_{\mb{B}}^{(k-1)} + \mb{Y}^{(k)\T}\left(\boldsymbol{\gamma}^{(k)\T}\odot\mb{A}^{(k)}\right),
\label{eq:GBk}
\end{eqnarray}
which implies the recursive computation of~(\ref{eq:Bkp1}) as 
\begin{eqnarray}
\mb{B}^{(k)} = \mb{G}_{\mb{B}}^{(k)}\left(\mb{V}_{\mb{B}}^{(k)}+\lambda\mb{D}_{2}^{(k-1)}\right)^{-1}. 
\label{eq:Bk}
\end{eqnarray}

The resulting algorithm, called \emph{Online BTD Reweighted Least Squares (RLS) (O-BTD-RLS)}, is tabulated as Algorithm~2.
\begin{table}
\centering
\title{Algorithm 2: The O-BTD-RLS algorithm}
 \begin{tabular}{|l|}
 \hline \\
  {\bf Algorithm 2}: The O-BTD-RLS algorithm\\ \hline 
  Input: $\bc{Y}$ in a streaming manner, $\xi,\lambda,\mu,R_{\mathrm{ini}},L_{\mathrm{ini}}$ \\
  Output: Best (in the sense of~(\ref{eq:ewminp})) BTD approximation of $\bc{Y}^{(k)}$ \\
  Initialize $\mb{A}^{(0)},\mb{B}^{(0)}, \mb{C}^{(0)},\mb{V}_{\mb{A}}^{(0)},\mb{V}_{\mb{B}}^{(0)},\mb{G}_{\mb{A}}^{(0)},\mb{G}_{\mb{B}}^{(0)}$ from Algorithm 1  \\
  \bf{for} $k=1,2,\ldots$\\
	  \hspace{0.3cm} Compute $\mb{D}_1^{(k-1)},\mb{D}_2^{(k-1)}$ from (\ref{eq:D1}) and (\ref{eq:D2}) \\
	  \hspace{0.3cm} Compute $\mb{S}^{(k-1)}$ from $\mb{A}^{(k-1)}$ and $\mb{B}^{(k-1)}$ (cf. (\ref{eq:X3})) \\
	   \hspace{0.3cm} $\boldsymbol{\gamma}^{(k)} \leftarrow \left(\mb{S}^{(k-1)\T}\mb{S}^{(k-1)}+\mu\mb{D}_{1}^{(k-1)}\right)^{-1}\mb{S}^{(k-1)\T}\mb{y}^{(k)}$
	   \\
	   \hspace{0.3cm} $\mb{V}_{\mb{A}}^{(k)} \leftarrow \xi\mb{V}_{\mb{A}}^{(k-1)} + \left(\mb{B}^{(k-1)}\odot\boldsymbol{\gamma}^{(k)\T}\right)^{\T}\left(\mb{B}^{(k-1)}\odot\boldsymbol{\gamma}^{(k)\T}\right)$
	   \\
	   \hspace{0.3cm} $\mb{G}_{\mb{A}}^{(k)} \leftarrow \xi\mb{G}_{\mb{A}}^{(k-1)} + \mb{Y}^{(k)}\left( \mb{B}^{(k-1)}\odot\boldsymbol{\gamma}^{(k)\T}\right)$
	   \\
    \hspace{0.3cm} $\mb{A}^{(k)} \leftarrow \mb{G}_{\mb{A}}^{(k)}\left(\mb{V}_{\mb{A}}^{(k)}+\lambda\mb{D}_{2}^{(k-1)}\right)^{-1}$ \\
    \hspace{0.3cm} $\mb{V}_{\mb{B}}^{(k)} \leftarrow \xi\mb{V}_{\mb{B}}^{(k-1)} + \left(\boldsymbol{\gamma}^{(k)\T}\odot\mb{A}^{(k)}\right)^{\T}\left(\boldsymbol{\gamma}^{(k)\T}\odot\mb{A}^{(k)}\right)$
	   \\
	   \hspace{0.3cm} $\mb{G}_{\mb{B}}^{(k)} \leftarrow \xi\mb{G}_{\mb{B}}^{(k-1)} + \mb{Y}^{(k)\T}\left(\boldsymbol{\gamma}^{(k)\T}\odot\mb{A}^{(k)}\right)$
	   \\
    \hspace{0.3cm} $\mb{B}^{(k)} \leftarrow \mb{G}_{\mb{B}}^{(k)}\left(\mb{V}_{\mb{B}}^{(k)}+\lambda\mb{D}_{2}^{(k-1)}\right)^{-1}$ \\
    \bf{end} \\
    \hline
 \end{tabular}
\end{table}
 The common in the OTF literature initialization practice, namely initializing the online scheme with the result of applying the batch algorithm in the first few slices, can also be employed here. 
%As a result, it exhibits very low sensitivity to initialization \cite{zvbjd16}, in view of the robustness of Algorithm 1 used to compute the initial BTD.
It should be stressed that, in the proposed algorithm, we respect the BTD model structure throughout, in contrast to works like~\cite{ns09} where the KR product structure of the slowly-varying part is only taken into account at the end of each recursion.

%Convergence: as in online dictionary learning~\cite{mbps10}; Kalman-like!

\subsection{Complexity analysis}
\label{sec:complexity}

Thanks to its recursive nature, O-BTD-RLS is much more efficient in terms of memory requirements than its batch counterpart. One can readily verify that, in the most practical case of big low-rank tensors, that is $I,J\gg R,L$, it requires a total storage of $\mathcal{O}([2(I+J)L+IJ]R)$ floating-point numbers, or roughly $\mathcal{O}(IJR)$. In contrast, the memory complexity of the BTD-IRLS algorithm for processing such a tensor with $I,J,K\gg R,L$ is of the order of $\mathcal{O}((IJ+(JK+KI)L)R)$ for storing the involved quantities, and one should add to this the $IJK$ places needed to store the entire tensor. As expected, the memory efficiency of the online w.r.t. the batch method increases with the size of the evolving mode, $K$.

In a manner analogous to that for BTD-HIRLS~\cite[Appendix~C]{rkg21}, it turns out that the per-iteration computational requirements of BTD-IRLS are of the order of $\mathcal{O}(IJKLR)$, again for the case of a big low-rank $I\times J\times K$ tensor. The computational cost of an O-BTD-RLS recursion is roughly estimated as follows. $\mb{D}_1$ and $D_2$ need $\mathrm{O}(R)$ and $\mathcal{O}((I+J)LR)$ multiplications and divisions, including those for the calculation of the square roots. The computation of $\mb{S}$ requires $IJLR$ multiplications and that of its Grammian, $\mb{S}^{\T}\mb{S}$ $(I+J+1)(LR)^2$ multiplications by virtue of \cite[Eq.~(25)]{rkg21}. For the update of $\boldsymbol{\gamma}$, $\mathcal{O}(IJR+R^3/3+R^2)$ multiplications are required. $(I+J)LR$ multiplications are needed to compute $\boldsymbol{\gamma}^{\T}\odot\mb{A}$ and $\mb{B}\odot\boldsymbol{\gamma}^{\T}$. In addition, we need $(2IJ+I+J)LR$ multiplications for $\mb{G}_{\mb{A}}$ and $\mb{G}_{\mb{B}}$. $(I+J+2)(LR)^2$ multiplications for $\mb{V}_{\mb{A}}$ and $\mb{V}_{\mb{B}}$, and $\mathcal{O}((I+J)(LR)^2+2(LR)^3/3)$ multiplications for computing $\mb{A}$ and $\mb{B}$. The total per-iteration cost is $\mathcal{O}(IJRL)$, that is, roughly the cost of analyzing the entire tensor in one batch divided by $K$. 

%time and memory complexity, compared to~\cite{gp20} and~\cite{rkg21}.

The complexity (in terms of memory and computation) of both schemes can of course be reduced if symmetries and/or sparsity (via efficiently computing the matricized tensor times Khatri-Rao products (MTTKRP)~\cite{sdfhpf17,gp20}) are also taken into account. 

%Approximate Grammians by diagonals; what is the accuracy loss?~\cite{k19}

\section{Simulation Results}
\label{sec:sims}

In this section, we evaluate the performance of the O-BTD-RLS algorithm in selecting and computing/tracking the correct BTD model for a given tensor. Comparisons with its batch counterpart are included, in terms of both approximation accuracy and time efficiency. We consider two experiments. The first, with a tensor of relatively large third mode which is described by a single BTD model in its entirety, corresponds to a big data scenario, where the tensor has to be analyzed incrementally due to its large size that prevents its storage and analysis as a whole. In this case, we can compare the batch and online schemes in terms of the error incurred when approximating the tensor and the time this requires. In the second experiment, a data-streaming scenario is considered, in which the tensor only becomes available on a slice-by-slice basis and its underlying BTD model may be time varying. This allows us to evaluate the tracking ability of the online algorithm. In both cases, we also assess the ability of the BTD-IRLS method to reveal the BTD ranks. O-BTD-RLS is initialized with the results of the batch scheme applied on the first few slices of the tensor. We show that the selected model is tracked by the online variant. 

We consider an $I\times J \times K$ tensor
$\bc{Y}=\bc{X}+\sigma\bc{N}$, where $\bc{X}$ is built as in~(\ref{eq:BTD}), $\bc{N}$ contains zero-mean, independent and identically distributed (i.i.d.) Gaussian entries of unit variance, and $\sigma$ is set so that we get a given signal-to-noise ratio (SNR), with SNR in~dB defined as $\mathrm{SNR}=10\log_{10}\|\bc{X}\|_{\F}^2/(\sigma^2\|\bc{N}\|_{\F}^2)$. The entries of $\mb{A},\mb{B},\mb{C}$ are also i.i.d. samples from the standard Gaussian distribution. 

\emph{Experiment~1:} We set $I=40$, $J=35$ and $K=1250$. The true $R$ is~5 and all $L_r$s are equal to~4. BTD-IRLS was run with overestimates of the true ranks, namely $R = 10$ and $L = 10$ for all blocks terms, using random initialization.\footnote{In the light of the robustness to initialization of BTD-HIRLS demonstrated in~\cite{rkg21}, we expect that BTD-IRLS will also enjoy such a desirable property. Further experimentation is required to confirm this.} O-BTD-RLS was initialized with the result of BTD-IRLS on the first~50 frontal slices and incrementally processed the remaining 1200~slices, with $\xi = 1$. For BTD-IRLS, $\lambda$ was selected as $L(I+J)\hat{\sigma}$, where $\hat{\sigma}$ is an estimate of the noise standard deviation (in our experiments taken equal to the true $\sigma$), and either $\mu = 0.75KR\hat{\sigma}$ (at SNR=5~dB) or $\mu = 2KR\hat{\sigma}$ (at SNR=10 and 15~dB). The regularization parameters of O-BTD-RLS were selected for simplicity equal to each other, $\lambda = \mu = L(I+J)\hat{\sigma}$.
Table~\ref{tab:noise_exp} shows the Relative Error (RE), $\left\|\bc{Y}-\bc{\hat{X}}\right\|_{\F}/\left\|\bc{Y}\right\|_{\F}$ and the normalized mean squared error (NMSE) over the blocks,  $\frac{1}{R}\sum^R_{r=1} \frac{\|\mathbf{A}_r\mathbf{B}^{\T}_r\circ \mathbf{c}_r - \hat{\mathbf{A}}_r\hat{\mathbf{B}}^{\T}_r\circ \hat{\mathbf{c}}_r\|_{\F}^2}{\|\mathbf{A}_r\mathbf{B}^{\T}_r\circ \mathbf{c}_r\|_{\F}^2}$, at different SNR values, for both the batch and incremental cases, where $\bc{\hat{X}}$ stands for the approximation of $\bc{Y}$ given by the computed BTD model and $(\mathbf{A,B,C}),$  $(\hat{\mathbf{A}},\hat{\mathbf{B}},\hat{\mathbf{C}})$ denote the true and the estimated BTD factors, respectively. For computing the NMSE, the Hungarian algorithm was employed to resolve permutation ambiguities (as in~\cite{rkg21}). In each case, the median over 50~independent realizations of the experiment is given.  
\begin{table*}
    \centering
        \caption{Comparison of BTD-IRLS and O-BTD-RLS in Terms of Relative Error, NMSE over the Blocks and Average Run-time at Different SNR Values.}
    \begin{tabular}{l|c|c|c|c|c|c|c|c|} 
    \cline{2-8} 
   & \multicolumn{6}{|c|}{SNR (dB)}  & \multirow{3}{*}{Average
   run-time (s)} \\ \cline{2-7}
     & \multicolumn{2}{c|}{5} & \multicolumn{2}{c|}{10} & \multicolumn{2}{c|}{15}  & \\ \cline{2-7} 
     & NSE ($\times10^{-3}$) & RE & NMSE ($\times10^{-3}$) & RE & NMSE ($\times10^{-3}$) & RE & \\ \hline
      \multicolumn{1}{|c|}{BTD-IRLS}   & $0.5$ & $0.2703$ & $0.15$   & $0.1560$ & $0.05$ &  $0.085$ &10.8\\ \hline
   \multicolumn{1}{|c|}{O-BTD-RLS}  & $1.2$ & $0.2945$ & $0.4$   & $0.16$ & $0.2$ & $0.0886$ &$2.8$ \\ \hline 
    \end{tabular}
    \label{tab:noise_exp}
\end{table*}
Clearly, the online scheme achieves an approximation accuracy close to that of the batch one (especially in terms of the RE), and the performance gap is diminished as SNR increases. %\textcolor{red}{However, NMSE values are not close. That of the onine is almost 4 times the one of the batch method. We must make sure that the model is correctly identified (at least at very high SNR). In addition, we do not evaluate the rank-revealing ability of the batch scheme, as we did before. To this end, many more realizations are needed, with both different X and N, again as we have previously done.} 
In addition, this is achieved much more efficiently in terms of the average run-time, with O-BTD-RLS requiring only {2.8~sec} to process the whole tensor, as compared to {the~10.8~sec} needed by BTD-IRLS. All experiments were run in a MacBook Pro, 2.6 GHz 6-Core Intel Core i7, 16 GB 2667 MHz DDR4 using Matlab v2019b.

\begin{figure}
\centerline{\includegraphics[width=0.5\textwidth]{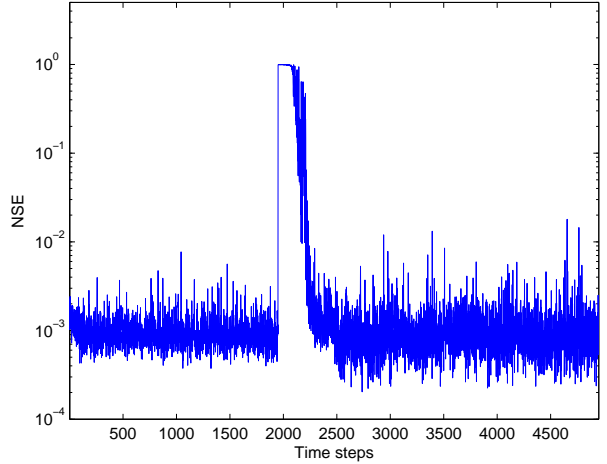}}
\caption{Normalized squared error per frontal slice  of O-BTD-RLS vs. time steps, at SNR=10~dB. The model changes abruptly at the $k=2001$ step.}
\label{fig:nse_tv}
\end{figure}

\emph{Experiment~2:} This experiment aims at assessing the ability of O-BTD-RLS to track changes of the BTD model. To this end, we consider a $40\times 35 \times 5000$ tensor resulting by concatenating two smaller tensors with dimensions $40\times 35 \times 2000$ and $40\times 35 \times 3000$. The true ranks $R$ and $L_r$s of the first tensor are~5 and~4 and those of the second one are~4 and~2. This implies an  abrupt change of the underlying BTD model at the $k=2001$ step. Noise is added for an SNR of~10~dB. Since the model is now time varying, a fading memory effect is simulated by setting $\xi$ to~0.985 in O-BTD-RLS. All other specifications are as in the previous experiment. The normalized squared error (NSE) per frontal slice of one run of the online algorithm is plotted in Fig.~\ref{fig:nse_tv} as a function of the number of update steps, namely $\mathrm{NSE}(k)=\left\|\bc{X}(:,:,k) - \mb{A}^{(k)}(\mathrm{diag}(\boldsymbol{\gamma}^{(k)})\otimes \mb{I}_L)\mb{B}^{(k)\T}\right\|_{\F}^{2}/\left\|\bc{X}(:,:,k)\right\|_{\F}^{2}$. 
Note that the algorithm immediately recognizes the model change and, after a few steps required to collect information about the new model, re-adapts fast back to the new BTD model.

\section{Conclusion}
\label{sec:concls}

The problem of rank-$(L_r,L_r,1)$ BTD model selection and tracking was studied in this work \emph{for the first time}, based on the idea of imposing column sparsity jointly on the factors and estimating the ranks as the numbers of factor columns of nonnegligible magnitude. An online method of the alternating reweighted least squares (RLS) type was developed, on the basis of a newly introduced rank-revealing batch scheme. It was shown to be computationally efficient and fast converging, with a modeling capacity comparable to that of its batch counterpart and the ability to track models that change abruptly in time. The proposed online scheme was shown to be memory efficient while its time efficiency was demonstrated to be considerably higher than that of its batch counterpart. The effectiveness in both selecting and tracking the correct BTD model was clearly demonstrated via simulation results. 

Future research can be directed towards extending the novel online BTD algorithm to incorporate constraints and side information~\cite{nmgt18} and perform completion~\cite{rgk21} and DL~\cite{zeb18} tasks. Modifications necessary for solving large-scale problems (using, for example, sampling/sketching)~\cite{gpp18,sgltu20} or being robust (to outliers)~\cite{rgk20,cdpm20} are also worth exploring. Coupled~\cite{ktwc19,wwzzz21} and Bayesian~\cite{grk21} versions should also be developed, greatly widening its applications spectrum.
The convergence analysis of the proposed online algorithm is also left as a future work. Note that, as pointed out in the introduction, sich a task is far from being easy for this kind of methods. It is expected that ideas recently employed to develop provable online CPD schemes (e.g., \cite{sln21}) and knowledge available for online sparsity learning~\cite[Chapter~10]{t20} will be most helpful in this direction. The effect of having fixed vs. time-varying regularization parameters is also worth investigating (cf.~\cite{dsbm20}).

\section{Acknowledgments}
The work of E. Kofidis has been partly supported by the University of Piraeus Research Center. P. V. Giampouras is supported by the European Union under the Horizon 2020 Marie Sk{\l}odowska-Curie Global Fellowship program:  HyPPOCRATES-H2020-MSCA-IF-2018, Grant Agreement Number: 844290.

\appendices

\section{Derivation of Eqs.~(\ref{eq:Ck+1}), (\ref{eq:Akp1}), (\ref{eq:Bkp1})}
\label{sec:proofs}

Call $f_{\boldsymbol{\gamma}}(\boldsymbol{\gamma})$ the objective function of~(\ref{eq:objCkp1}), i.e., 
\begin{align}
f_{\boldsymbol{\gamma}}(\boldsymbol{\gamma}) =& \frac{1}{2}\left\|\mb{y}^{(k)}-\mb{S}^{(k-1)}\boldsymbol{\gamma}\right\|_{2}^2 + \nonumber \\
&\mu\sum_{r=1}^{R}\sqrt{\xi\left\|\mb{\Xi}^{(k-1)\frac{1}{2}}\mb{c}_{r}^{(k-1)}\right\|_{2}^2 + \gamma_r^2 + \eta^2}.
\label{eq:objc}    
\end{align}
To minimize this we follow a MM approach~\cite{hl04}. Namely, we instead minimize a surrogate function, which is selected here to be a second-order Taylor approximation of $f_{\boldsymbol{\gamma}}(\boldsymbol{\gamma})$ around $\boldsymbol{\gamma}=\mb{0}$:
\begin{align}
g_{\boldsymbol{\gamma}}(\boldsymbol{\gamma}) = f_{\boldsymbol{\gamma}}(\mb{0}) + \boldsymbol{\gamma}^{\T}\nabla_{\boldsymbol{\gamma}}f_{\boldsymbol{\gamma}}(\mb{0})+\frac{1}{2}\boldsymbol{\gamma}^{\T}\tilde{\mb{H}}_{\boldsymbol{\gamma}}^{(k-1)}\boldsymbol{\gamma},
\label{eq:surc}    
\end{align}
with $\nabla_{\boldsymbol{\gamma}}f_{\boldsymbol{\gamma}}(\mb{0})=-\mb{S}^{(k-1)\T}\mb{y}^{(k)}$ and $\tilde{\mb{H}}_{\boldsymbol{\gamma}}^{(k-1)}$ being an approximation of the Hessian matrix of $f_{\boldsymbol{\gamma}}$ at $\mb{0}$, $\mb{H}_{\boldsymbol{\gamma}}^{(k-1)}$:
\begin{eqnarray}
\tilde{\mb{H}}_{\boldsymbol{\gamma}}^{(k-1)} = \mb{S}^{(k-1)\T}\mb{S}^{(k-1)} + \mu\mb{D}_1^{(k-1)}. 
\label{eq:apHc}
\end{eqnarray}
The reason why we consider the Taylor approximation around $\mb{0}$ and not $\boldsymbol{\gamma}^{(k-1)}$ is that the variable $\boldsymbol{\gamma}^{(k)}$ ($k$th row of $\mb{C}$) is not related in any way to $\boldsymbol{\gamma}^{(k-1)}$ ($(k-1)$st row of $\mb{C}$).

Note from~(\ref{eq:surc}) that $g_{\boldsymbol{\gamma}}(\mb{0})=f_{\boldsymbol{\gamma}}(\mb{0})$ and $\tilde{\mb{H}}_{\boldsymbol{\gamma}}^{(k-1)}$ is positive definite. In addition, as we will show below, the matrix  $\tilde{\mb{H}}_{\boldsymbol{\gamma}}^{(k-1)}-\mb{H}_{\boldsymbol{\gamma}}^{(k-1)}$ is positive semi-definite. These properties ensure that  $g_{\boldsymbol{\gamma}}(\boldsymbol{\gamma})$ is a majorizing function of $f_{\boldsymbol{\gamma}}(\boldsymbol{\gamma})$~\cite{hrlp16}. Indeed, it is not difficult to show that the Hessian of $f_{\boldsymbol{\gamma}}(\boldsymbol{\gamma})$ is given by 
\begin{eqnarray}
\mb{H}_{\boldsymbol{\gamma}}^{(k-1)} = \mb{S}^{(k-1)\T}\mb{S}^{(k-1)} + \mu\mb{D}_{\boldsymbol{\gamma}}^{(k-1)}, 
\label{eq:Hc}
\end{eqnarray}
where the entries of the $R\times R$ diagonal matrix $\mb{D}_{\boldsymbol{\gamma}}^{(k-1)}$ are expressed as 
\begin{eqnarray}
\mb{D}_{\boldsymbol{\gamma}}^{(k-1)}(r,r) = \left(\xi\left\|\mb{\Xi}^{(k-1)\frac{1}{2}}\mb{c}_{r}^{(k-1)}\right\|_{2}^2 + \gamma_r^2 + \eta^2\right)^{-\frac{1}{2}} \nonumber \\
- \gamma_r^2\left(\xi\left\|\mb{\Xi}^{(k-1)\frac{1}{2}}\mb{c}_{r}^{(k-1)}\right\|_{2}^2 + \gamma_r^2 + \eta^2\right)^{-\frac{3}{2}}.
\label{eq:Dc}
\end{eqnarray}
From~(\ref{eq:apHc}), (\ref{eq:D1}), (\ref{eq:Hc}) and~(\ref{eq:Dc}), we easily deduce that $\tilde{\mb{D}}_{\boldsymbol{\gamma}}^{(k-1)} \triangleq\tilde{\mb{H}}_{\boldsymbol{\gamma}}^{(k-1)}-\mb{H}_{\boldsymbol{\gamma}}^{(k-1)}$ is a diagonal matrix with non-negative diagonal elements given by  
\begin{eqnarray}
\tilde{\mb{D}}_{\boldsymbol{\gamma}}^{(k-1)}(r,r) = \mu\left(\xi\left\|\mb{\Xi}^{(k-1)\frac{1}{2}}\mb{c}_{r}^{(k-1)}\right\|_{2}^2 + \eta^2\right)^{-\frac{1}{2}} \nonumber \\
-\mu\left(\xi\left\|\mb{\Xi}^{(k-1)\frac{1}{2}}\mb{c}_{r}^{(k-1)}\right\|_{2}^2 + \gamma_r^2 + \eta^2\right)^{-\frac{1}{2}} \nonumber \\
+ \gamma_r^2\mu\left(\xi\left\|\mb{\Xi}^{(k-1)\frac{1}{2}}\mb{c}_{r}^{(k-1)}\right\|_{2}^2 + \gamma_r^2 + \eta^2\right)^{-\frac{3}{2}}.
\label{eq:tDc}
\end{eqnarray}  
Note that the $\tilde{\mb{D}}_{\boldsymbol{\gamma}}^{(k-1)}(r,r)$'s are all positive for $\boldsymbol{\gamma}\neq \mb{0}$ and become zero at $\boldsymbol{\gamma}=\mb{0}$, i.e., at the point of approximation. 

Besides majorizing the objective function $f_{\boldsymbol{\gamma}}$, the approximation function $g_{\boldsymbol{\gamma}}$ has the same first-order behavior with $f_{\boldsymbol{\gamma}}$ at $\mb{0}$~\cite{hrlp16}. As a result, the vector $\boldsymbol{\gamma}$ that minimizes $g_{\boldsymbol{\gamma}}$ will at least guarantee some descent of the original objective $f_{\boldsymbol{\gamma}}$. It can be easily shown that the minimization of $g_{\boldsymbol{\gamma}}$ leads to~(\ref{eq:Ck+1}).

Consider now the objective function for the $\mb{A}$ sub-problem, 
\begin{align}
f_{\mb{A}}(\mb{A}) =&\frac{1}{2}\left\|{}_{J}\mb{\Xi}^{(k)\frac{1}{2}}\left(\mb{Y}_{(1)}^{(k)\T}-\mb{P}^{(k)}\mb{A}^{\T}\right)\right\|_{\F}^2 + \nonumber \\
&\lambda\sum_{r=1}^{R}\sum_{l=1}^{L}\sqrt{\|\mb{a}_{r,l}\|_{2}^2 + \|\mb{b}_{r,l}^{(k-1)}\|_{2}^2 + \eta^2},
\label{eq:objfA}    
\end{align}
and define a surrogate function, $g_{\mb{A}}(\mb{A})$, as a second order Taylor approximation of $f_{\mb{A}}(\mb{A})$ around $\mb{A}^{(k-1)}$, the estimate of $\mb{A}$ computed at step $k-1$,  
\begin{align}
g_{\mb{A}}(\mb{A}) &= f_{\mb{A}}(\mb{A}^{(k-1)}) + \text{tr}\{(\mb{A}-\mb{A}^{(k-1)})  
\nabla_{\mb{A}}f_{\mb{A}}(\mb{A}^{(k-1)})\} \nonumber \\ 
&+ \frac{1}{2}\text{vec}(\mb{A}-\mb{A}^{(k-1)})^{\T}\tilde{\mb{H}}_{\mb{A}^{(k-1)}}\text{vec}(\mb{A}-\mb{A}^{(k-1)}),
\label{eq:surA}
\end{align}
where $\tilde{\mb{H}}_{\mb{A}^{(k-1)}}$ is an approximation of the $IRL \times IRL$  Hessian $\mb{H}_{\mb{A}}$ of $f_{\mb{A}}$ at $\mb{A}^{(k-1)}$, defined as  
\begin{eqnarray}
\tilde{\mb{H}}_{\mb{A}^{(k-1)}} = \mb{P}^{(k)\T}{}_{J}\mb{\Xi}^{(k)}\mb{P}^{(k)}+\lambda\mb{D}_{2}^{(k-1)}. 
\label{eq:apHA}
\end{eqnarray}
It is now clear that $g_{\mb{A}}(\mb{A}^{(k-1)}) = f_{\mb{A}}(\mb{A}^{(k-1)})$ and $\tilde{\mb{H}}_{\mb{A}^{(k-1)}}$ is positive definite. In addition, working as in~\cite[Appendix~B]{grk19a}, it can be shown that $\tilde{\mb{H}}_{\mb{A}^{(k-1)}}-\mb{H}_{\mb{A}^{(k-1)}}$ is also positive definite, which verifies that $g_{\mb{A}}$ majorizes  $f_{\mb{A}}$ around $\mb{A}^{(k-1)}$. It is then straightforward to show that the unique minimizer of the quadratic function $g_{\mb{A}}$ is given by~(\ref{eq:Akp1}). 

In a similar way, we can arrive at the update equation~(\ref{eq:Bkp1}) for the $\mb{B}$ factor.  

\bibliographystyle{IEEEtran}
\bibliography{IEEEabrv,refs}

\end{document}